\documentclass[leqno, notitlepage]{amsart}
\usepackage{amssymb}

\DeclareMathOperator{\End}{End}
\DeclareMathOperator{\Hom}{Hom}
\DeclareMathOperator{\res}{res}
\DeclareMathOperator{\im}{im}

\newcommand{\sss}{\scriptscriptstyle}

\swapnumbers
\theoremstyle{plain} 
\newtheorem*{theorem}{Theorem}
\newtheorem*{proposition}{Proposition}

\theoremstyle{remark}
\newtheorem*{remark}{Remark}

\numberwithin{equation}{subsection}

\numberwithin{enumi}{subsection}

\begin{document}

\title{A generalization of Hall polynomials to
ADE case}

\author{ Anton Malkin }
\address{ Department of Mathematics, Yale University,
P.O. Box 208283, New Haven, CT 06520-8283 }
\email{anton.malkin@yale.edu}

\begin{abstract}
Certain computable 
polynomials are described whose leading
coefficients are equal to multiplicities
in the tensor product decomposition for 
representations of a Lie algebra of
ADE type. 
\end{abstract}

\maketitle

\subsection{Hall theorem}
\label{HallTheorem}

Given a prime number $p$ and a positive 
integer $l$ let $\mathbb{F}_p$ be a field 
with $p$ elements, $\overline{\mathbb{F}}_p$ be
an algebraic closure of $\mathbb{F}_p$, 
and $\mathbb{F}_{p^l}$ be the subfield of
$\overline{\mathbb{F}}_p$ consisting of all 
$x \in \overline{\mathbb{F}}_p$  
such that $x^{p^l} = x$. Thus 
$\mathbb{F}_{p^l}$ is a field with
$p^l$ elements.

Let $k$ be a field, and let $k [[t]]$
be the ring of formal power series with 
coefficients in $k$ considered as a $k$-algebra.
To a representation $\pi$ of
$k [[t]]$ in a finite dimensional
$k$-linear space $D$
one can associate a partition $\lambda$
as follows: $\lambda_i = 
\dim_k \bigl( \im \pi (t^{i-1}) /
\im \pi (t^i) \bigr)$. This partition
is called the \emph{type} of the representation
$\pi$ (in \cite{Macdonald1995} the conjugate
partition is called the type of $\pi$).
Two representations are isomorphic if and 
only if they are of the same type.

Given a partition $\lambda$ let 
$\lambda_1 \geq \lambda_2 \geq \ldots$
be the parts of $\lambda$ and let
$| \lambda | = \sum_i \lambda_i$.
Let $\mathfrak{O}_{\lambda} (k)$ be the
set of all representations of $k [[t]]$
in $k^{| \lambda |}$ of type $\lambda$.
If $k$ is algebraically closed then
$\mathfrak{O}_{\lambda} (k)$ is an
irreducible smooth quasi-projective
variety of dimension $o_{\lambda} =
\sum_{i \neq j} \lambda_i \lambda_j$.
Let $\mathbf{o}_{\lambda} (p^l)$ be the cardinal
of $\mathfrak{O}_{\lambda} (\mathbb{F}_{p^l})$.
It is known (cf. \cite[Chapter II]{Macdonald1995})
that $\mathbf{o}_{\lambda} (p^l)$
is a monic polynomial function of $p^l$ 
of degree $o_{\lambda}$.

Let $\lambda$, $\mu^1$, \ldots , $\mu^n$, 
be partitions such that
$| \lambda | = | \mu^1 | + \ldots + | \mu^n |$.
Let 
$\mathfrak{N}_{\mu^1 \ldots \mu^n}^{\lambda} 
(k)$ be the set of  
all pairs $(\pi , \underline{D})$ consisting 
of a representation $\pi$ of 
$k [[t]]$ in $k^{| \lambda |}$ 
of type $\lambda$
and an $n$-step filtration
$\underline{D} = 
( \{ 0 \} = D^0 \subset D^1 \subset
\ldots \subset D^n = k^{| \lambda |})$ 
of $k^{| \lambda |}$ 
by subrepresentations of $\pi$ such that
$\pi |_{\sss D^a / D^{a-1}}$
has type $\mu^a$ for $a = 1, \ldots , n$. 
If $k$ is algebraically closed then 
$\mathfrak{N}_{\mu^1 \ldots \mu^n}^{\lambda} (k)$
is a quasi-projective variety over $k$. The set
$\mathfrak{N}_{\mu^1 \ldots \mu^n}^{\lambda} 
(\mathbb{F}_{p^l})$ can be identified
with the set of $\mathbb{F}_{p^l}$-rational
points in the variety 
$\mathfrak{N}_{\mu^1 \ldots \mu^n}^{\lambda} 
(\overline{\mathbb{F}}_p)$.

Steinitz \cite{Steinitz1901} and 
Hall \cite{Hall1959} observed that the sets 
$\mathfrak{N}_{\mu^1 \ldots \mu^n}^{\lambda} (k)$
have some relation to the theory of symmetric
functions, and thus to the representation
theory of the general linear and symmetric 
groups. Given a partition $\lambda$ such that
the number of non-zero parts of $\lambda$ is
less than $(N+1)$ let $L_{\lambda}$ 
denote the irreducible
representation of $GL (N, \mathbb{C})$
with the highest weight $\lambda$, and
given an arbitrary partition $\lambda$ let
$\rho_{\lambda}$ be the irreducible  
representation of $S_{| \lambda |}$ 
over $\mathbb{C}$ associated to $\lambda$.
Let $c_{\mu^1 \ldots \mu^n}^{\lambda}$ be
a non-negative integer defined in either of
the following two 
(equivalent) ways:
\begin{align}\label{HallTensor}
c_{\mu^1 \ldots \mu^n}^{\lambda} &=
\dim_{\sss \mathbb{C}} 
\Hom_{\sss GL (N, \mathbb{C})} 
(L_{\lambda}, L_{\mu^1} \otimes
\ldots \otimes L_{\mu^n}) \; , \;
\text{ or }
\\\label{HallRestriction}
c_{\mu^1 \ldots \mu^n}^{\lambda} &=
\dim_{\sss \mathbb{C}}
\Hom_{\sss S_{| \mu^1 |} \times \ldots \times
S_{| \mu^n |}} 
(\res^{S_{| \mu^1 |} \times \ldots \times
S_{| \mu^n |}}_{S_{| \lambda |}}
\rho_{\lambda}, \rho_{\mu^1} \otimes
\ldots \otimes \rho_{\mu^n}) \; .
\end{align}
In \eqref{HallTensor} it is assumed
that the number of non-zero parts of each of
the partitions $\lambda ,
\mu^1 , \ldots , \mu^n ,$ is
less that $(N+1)$. 

Let $v^1 , \ldots , v^n 
\in \mathbb{Z}_{\geq 0}$,
$u = \sum_{a=1}^n v^a$.
Assume $k$ is algebraically closed.
Let $f_{v^1 \ldots v^n}^{u}=
\sum_{a < b} v^a v^b$
be the dimension of the projective 
variety of all $n$-step partial flags in 
$k^u$ with dimensions of the 
subfactors given by
$v^1 , \ldots , v^n$.

\begin{theorem}\noindent\par
\begin{enumerate}
\item\label{HallDim}
If $k$ is algebraically closed
then the variety 
$\mathfrak{N}_{\mu^1 \ldots \mu^n}^{\lambda} (k)$
is empty or has pure dimension
\begin{equation}\nonumber
\dim 
\mathfrak{N}_{\mu^1 \ldots \mu^n}^{\lambda} (k) =
f_{| \mu^1 | \ldots | \mu^n |}^{| \lambda |} +
\tfrac{1}{2} \bigl(
o_{\lambda} +
o_{\mu^1} + \ldots + o_{\mu^n} 
\bigr) \; .
\end{equation}
\item\label{HallIrr}
If $k$ is algebraically closed
then the number of irreducible components of 
$\mathfrak{N}_{\mu^1 
\ldots \mu^n}^{\lambda} (k)$
is equal to $c_{\mu^1 \ldots \mu^n}^{\lambda}$.
\item\label{HallPolynomial}
There exists a computable polynomial 
$\mathbf{n}_{\mu^1 \ldots \mu^n}^{\lambda}$
with integer coefficients
such that the cardinal of the set 
$\mathfrak{N}_{\mu^1 \ldots \mu^n}^{\lambda} 
(\mathbb{F}_{p^l})$
is equal to 
$\mathbf{n}_{\mu^1 
\ldots \mu^n}^{\lambda} (p^l)$
for any prime number $p$ and any 
positive integer $l$.
\item\label{HallDegree}
The polynomial 
$\mathbf{n}_{\mu^1 \ldots \mu^n}^{\lambda}$
is identically equal to zero or has degree
\begin{equation}\nonumber
\deg 
\mathbf{n}_{\mu^1 \ldots \mu^n}^{\lambda} =
f_{| \mu^1 | \ldots | \mu^n |}^{| \lambda| } +
\tfrac{1}{2} \bigl(
o_{\lambda} +
o_{\mu^1} + \ldots + o_{\mu^n} 
\bigr) \; .
\end{equation}
\item\label{HallLR}
The leading coefficient of
the polynomial
$\mathbf{n}_{\mu^1 \ldots \mu^n}^{\lambda}$ 
is equal to
$c_{\mu^1 \ldots \mu^n}^{\lambda}$.
\end{enumerate}
\end{theorem}

\begin{remark}
$\mathbf{n}_{\mu^1 \ldots \mu^n}^{\lambda} 
(x) = \mathbf{o}_{\lambda} (x) \,
\mathbf{h}_{\mu^1 \ldots \mu^n}^{\lambda} 
(x)$, where 
$\mathbf{h}_{\mu^1 \ldots \mu^n}^{\lambda}$
is the Hall polynomial
(cf. \cite[Chapter II]{Macdonald1995}).
\end{remark}

\begin{proof}
Statement \ref{HallDim} is proven by
Spaltenstein \cite[II.5]{Spaltenstein1982};
\ref{HallPolynomial}, \ref{HallDegree},
and \ref{HallLR} are due to
Steinitz \cite{Steinitz1901} and
Hall \cite{Hall1959} (see also
\cite[Chapter II]{Macdonald1995});
\ref{HallIrr} follows from
\ref{HallDim} and \ref{HallLR}.
\end{proof}

Definitions \eqref{HallTensor} and
\eqref{HallRestriction} of 
$c_{\mu^1 \ldots \mu^n}^{\lambda}$ are
equivalent due to Schur-Weyl duality.
However if one wants to generalize 
Theorem \ref{HallTheorem} to 
reductive groups other than
$GL (N)$ there are various directions. 
For example certain subvarieties of
parabolic flag varieties play the role of 
$\mathfrak{N}^{\lambda}_{\mu^1 \ldots \mu^n}$
for restriction multiplicities
in the representation theory of Weyl groups.
In this note some other varieties 
are considered which are relevant to  
tensor product multiplicities. 
They
were independently described by Nakajima 
\cite{Nakajima2001p},
Varagnolo and Vasserot
\cite{VaragnoloVasserot2001},
and the author \cite{Malkin2000b}.

\subsection{A Lie algebra $\mathfrak{g}'$}
\label{LieAlgebraGPrime}

Let $\mathfrak{g}$ be a simple simply
laced Lie algebra over $\mathbb{C}$ and
$\mathfrak{t}$ be
a Cartan subalgebra of $\mathfrak{g}$.
Let $\mathfrak{g}' = 
\mathfrak{g} \oplus \mathfrak{t}$
(a reductive Lie algebra) and let
$\mathfrak{t}' = 
\mathfrak{t} \oplus \mathfrak{t} 
\subset
\mathfrak{g} \oplus \mathfrak{t}$
be a Cartan subalgebra of $\mathfrak{g}'$.
Let $I$ be the set of vertices
of the Dynkin graph of $\mathfrak{g}$.
Identify the weight lattice 
$\mathcal{Q}_{\mathfrak{g}}$
of $\mathfrak{g}$ with
$\mathbb{Z} [I]$ in such a way
that $i \in I \subset \mathbb{Z} [I]$
is the simple weight corresponding 
to the vertex $i$. Let
$\mathcal{Q}_{\mathfrak{g}}^{++} =
\mathbb{Z}_{\geq 0} [I]$ be the
set of positive weights. 
Fix a lattice $\mathcal{Q}_{\mathfrak{t}}$
in $\mathfrak{t}^{\ast}$ and identify
$\mathcal{Q}_{\mathfrak{t}}$ 
with $\mathbb{Z} [I]$ in some way.
Let $\mathcal{Q}^{++}_{\mathfrak{t}} =
\mathbb{Z}_{\geq 0} [I]$.
In what follows representations of
$\mathfrak{g}'$ are assumed to have 
the action of $\{ 0 \} \oplus \mathfrak{t}$ 
given by elements of the lattice 
$\mathcal{Q}_{\mathfrak{t}}$.
One has an injective map 
$\varkappa: 
\mathbb{Z} [I] \oplus \mathbb{Z} [I]
\rightarrow 
\mathcal{Q}_{\mathfrak{g}} \oplus
\mathcal{Q}_{\mathfrak{t}}$ given by 
$\varkappa ((\mathbf{u},\mathbf{v})) = 
(\mathbf{u}-\mathbf{v} , 
\mathbf{u}-\mathbf{v}+A \mathbf{v})$, 
where $A$ is the Cartan matrix of
$\mathfrak{g}$. Let
$\mathcal{Q}_{\mathfrak{g}'}$ be
the image of the map $\varkappa$.
In what follows 
$\mathcal{Q}_{\mathfrak{g}'}$
is called the \emph{weight lattice of
$\mathfrak{g}'$}.
For $\xi \in \mathcal{Q}_{\mathfrak{g}'}
\subset \mathcal{Q}_{\mathfrak{g}} \oplus
\mathcal{Q}_{\mathfrak{t}}$ let
$| \xi |$ denote the 
$\mathcal{Q}_{\mathfrak{t}}$-component
of $\xi$.
A weight $\xi \in 
\mathcal{Q}_{\mathfrak{g}'}$
is called \emph{integrable} if 
$\varkappa^{-1} ( \xi )
\in \mathbb{Z}_{\geq 0} [I] \oplus
\mathbb{Z}_{\geq 0} [I]
\subset \mathbb{Z} [I] \oplus
\mathbb{Z} [I]$. Let
$\mathcal{Q}^{+}_{\mathfrak{g}'}$ be
the set of integrable weights, and
let $\mathcal{Q}^{++}_{\mathfrak{g}'}=
\mathcal{Q}^{+}_{\mathfrak{g}'} \cap
(\mathcal{Q}^{++}_{\mathfrak{g}} \oplus
\mathcal{Q}^{++}_{\mathfrak{t}})$.
Elements of 
$\mathcal{Q}^{++}_{\mathfrak{g}'}$
are called
\emph{positive integrable}
weights. 

Given 
$\lambda \in \mathcal{Q}^{++}_{\mathfrak{g}'}$
let $L_{\lambda}$ be
an irreducible highest weight representation of
$\mathfrak{g}'$ with the highest weight 
$\lambda$.
The representation $L_{\lambda}$ 
is finite dimensional and all its
weights belong to
$\mathcal{Q}^{+}_{\mathfrak{g}'}$.
A finite dimensional representation of 
$\mathfrak{g}'$ is called integrable 
if all its irreducible components have
positive integrable highest weights. 
The category of
finite dimensional integrable representations
is closed with respect to tensor product.
Let $d_{\eta^1 \ldots \eta^n}^{\xi}$ be a
tensor product multiplicity:
\begin{equation}\nonumber
d_{\eta^1 \ldots \eta^n}^{\xi} =
\dim_{\sss \mathbb{C}}
\Hom_{\mathfrak{g}'} 
(L_{\xi}, L_{\eta^1} \otimes
\ldots \otimes L_{\eta^n}) \; ,
\end{equation}
where $\xi, \eta^1 , \ldots , \eta^n \in
\mathcal{Q}^{++}_{\mathfrak{g}'}$ and
$| \xi | = | \eta^1 | + \ldots + | \eta^n |$. 

\subsection{An algebra $\Tilde{\mathcal{F}}$}

Recall that $I$ is the set of vertices of the
Dynkin graph of $\mathfrak{g}$. Let
$H$ be the set of pairs consisting of
an edge of the Dynkin graph together with an
orientation of the edge. Given a field
$k$ let $\mathcal{F} (k)$ be the
path algebra of the Dynkin graph over $k$,
and let $\cdot$ denote 
the multiplication in 
$\mathcal{F} (k)$. 
Fix a map 
$\varepsilon : H \rightarrow \{ \pm 1 \}$
such that $\varepsilon (h) = - 
\varepsilon (\overline{h})$, where 
$\overline{h}$ denotes the same edge 
as $h$ but with the opposite orientation.
Let $\vartheta = \sum_{h \in H} 
\varepsilon (h) h \cdot \overline{h}
\in \mathcal{F} (k)$, and
let $\Tilde{\mathcal{F}} (k)$ be an 
associative $k$-algebra defined as follows.
As a $k$-linear space 
$\Tilde{\mathcal{F}} (k) =
\mathcal{F} (k) \oplus
\oplus_{i \in I} 
k u_i$, where $u_i$ are some symbols,
and the multiplication
$\circ$ in $\Tilde{\mathcal{F}} (k)$
is given by 
\begin{align}\nonumber
f &\circ f' = f \cdot \vartheta \cdot f' \; ,
&
f &\circ u_i = f \cdot [i] \; ,
\\\nonumber
u_i & \circ u_j = \delta_{ij} \, u_i \; ,
&
u_i &\circ f \; = [i] \cdot f \; ,
\end{align}
where $f , f' \in \mathcal{F} (k)$,
$i , j \in I$,
$[i] \in \mathcal{F} (k)$ is the
path of length zero beginning and ending at
a vertex $i \in I$,
and $\delta_{ij}$ is the Kronecker symbol.
The algebra $\Tilde{\mathcal{F}} (k)$ was 
introduced by Lusztig.
It is finitely generated as a $k$-algebra
(cf. \cite[2.2]{Lusztig2000a}).

Let $k^I$ be a semisimple 
$k$-algebra isomorphic to the
direct sum of $|I|$ copies 
of $k$. In what follows
$\mathcal{F} (k)$ and 
$\Tilde{\mathcal{F}} (k)$
are considered as $k^I$-algebras
with the embedding of the set of 
idempotents of $k^I$ into 
$\mathcal{F} (k)$ 
(resp. $\Tilde{\mathcal{F}} (k)$)
given by
$\{ [i] \}$ (resp. $\{ u_i \}$).

Let $\mathbf{D}$ be a $k^I$-module 
(i.e. $\mathbb{Z} [I]$-graded
$k$-linear space). A representation
of $\Tilde{\mathcal{F}} (k)$ in
$\mathbf{D}$ is a $k^I$-algebra 
homomorphism $\Tilde{\mathcal{F}} (k)
\rightarrow \End_{k^I} \mathbf{D}$. 
Let $\Dot{\mathbf{D}} = \mathcal{F} (k) 
\otimes_{k^I} \mathbf{D}$ 
(cf. \cite{Lusztig1998, Lusztig2000a}).
$\Dot{\mathbf{D}}$ is naturally a left
$\mathcal{F} (k)$-module.
Given a representation $\pi$
of $\Tilde{\mathcal{F}} (k)$ in
$\mathbf{D}$ let $\varpi_{\pi} \in
\Hom_{k^I} (\Dot{\mathbf{D}}, \mathbf{D})$ 
be given by 
$\varpi_{\pi} (f \otimes d) = 
\pi (f) d$, and let $\mathcal{K}_{\pi}$  
be the maximal $\mathcal{F} (k)$-submodule
of $\Dot{\mathbf{D}}$ contained in the 
kernel of $\varpi_{\pi}$. 
Let $\xi \in \mathcal{Q}_{\mathfrak{g}'}$
be given by $\xi = \varkappa (
(\mathbf{d} + \mathbf{v} - A \mathbf{v} , 
\mathbf{v}))$, where 
$\mathbf{d} = \dim_{k^I} \mathbf{D}$, 
$\mathbf{v} = \dim_{k^I} 
(\Dot{\mathbf{D}}/\mathcal{K}_{\pi})$
(note that $| \xi | =\mathbf{d}$).
Such $\xi$ is called the \emph{type} of $\pi$.
Two representations
of the same type are not necessarily
isomorphic.
Given $\xi \in \mathcal{Q}_{\mathfrak{g}'}$
let $\mathfrak{X}_{\xi} (k)$ be the set
of all representations of
$\Tilde{\mathcal{F}} (k)$ in 
$k^{| \xi |}$ of type $\xi$.
It is known (cf. 
\cite{Nakajima1998},
\cite{Lusztig1998, Lusztig2000a},
\cite{CrawleyBoevey2000})
that $\mathfrak{X}_{\xi} (k)$ is non-empty
if and only if
$\xi \in \mathcal{Q}^{++}_{\mathfrak{g}'}$,
and that if $k$ is algebraically closed then
$\mathfrak{X}_{\xi} (k)$ is empty or is an
irreducible smooth quasi-projective
variety over $k$ of dimension
$x_{\xi} = \sum_{i \in I} 
\mathbf{v}_i 
(2 | \xi | - A \mathbf{v}_i )$,
where $(\mathbf{u} , \mathbf{v})= 
\varkappa^{-1} ( \xi )$.

\subsection{Generalized Hall theorem}
\label{MultTheorem}

Let $\xi, \eta^1, \ldots , 
\eta^n \in \mathcal{Q}^{++}_{\mathfrak{g}'}$, 
$| \xi | = |\eta^1 | + \ldots + | \eta^n |$,
and let 
$\mathfrak{P}_{\eta^1 \ldots \eta^n}^{\xi} 
(k)$ be the set
of all pairs $(\pi, \underline{\mathbf{D}})$ 
consisting of a
representation $\pi$ of $\Tilde{\mathcal{F}}$ 
in $k^{| \xi |}$ of type $\xi$
and an $n$-step filtration 
$\underline{\mathbf{D}} = 
( \{ 0 \} = \mathbf{D}^0 
\subset \mathbf{D}^1 \subset
\ldots \subset \mathbf{D}^n = k^{| \xi |} )$
of $k^{| \xi |}$ by subrepresentations of 
$\pi$ such that
$\pi|_{\sss \mathbf{D}^a /\mathbf{D}^{a-1}}$
has type $\eta^a$ for $a = 1, \ldots , n$. 
If $k$ is algebraically closed then 
$\mathfrak{P}_{\eta^1 \ldots \eta^n}^{\xi} 
(k)$ is a quasi-projective variety over $k$. 
The variety
$\mathfrak{P}_{\eta^1 \ldots \eta^n}^{\xi} 
(k)$ is a variant of a
``multiplicity variety'' of \cite{Malkin2000b}. 
It is also implicitly
contained in \cite{Nakajima2001p} and
\cite{VaragnoloVasserot2001}.
The set
$\mathfrak{P}_{\eta^1 \ldots \eta^n}^{\xi} 
(\mathbb{F}_{p^l})$ can be identified
with the set of $\mathbb{F}_{p^l}$-rational
points in the variety 
$\mathfrak{P}_{\eta^1 \ldots \eta^n}^{\xi} 
(\overline{\mathbb{F}}_p)$.

Let $\mathbf{v}^1 , \ldots , \mathbf{v}^n 
\in \mathbb{Z}_{\geq 0} [I]$,
$\mathbf{u} = \sum_{a=1}^n \mathbf{v}^a$.
Assume $k$ is algebraically closed.
Let $g_{\mathbf{v}^1
\ldots \mathbf{v}^n}^{\mathbf{u}} =
\sum_{\substack{a < b \\ i \in I}} 
\mathbf{v}^a_i \mathbf{v}^b_i$
be the dimension of the projective 
variety of all $n$-step filtrations
of $k^{\mathbf{u}}$ by $k^I$-submodules
with dimensions of the 
subfactors given by
$\mathbf{v}^1 , \ldots , \mathbf{v}^n$. 

\begin{theorem}\noindent\par
\begin{enumerate}
\item\label{MultDim}
If $k$ is algebraically closed
then the variety 
$\mathfrak{P}_{\eta^1 \ldots \eta^n}^{\xi} 
(k)$ is empty or has pure dimension
\begin{equation}\nonumber
\dim 
\mathfrak{P}_{\eta^1 \ldots \eta^n}^{\xi} 
(k) =
g_{| \eta^1 | \ldots | \eta^n |}^{| \xi |} +
\tfrac{1}{2} \bigl(
x_{\xi} +
x_{\eta^1} + \ldots + x_{\eta^n} 
\bigr) \; .
\end{equation}
\item\label{MultIrr}
If $k$ is algebraically closed
then the number of irreducible components of 
$\mathfrak{P}_{\eta^1 \ldots \eta^n}^{\xi} (k)$
is equal to $d_{\eta^1 \ldots \eta^n}^{\xi}$.
\item\label{MultPolynomial}
There exists a computable polynomial 
$\mathbf{p}_{\eta^1 \ldots \eta^n}^{\xi}$
with integer coefficients
such that the cardinal of the set 
$\mathfrak{P}_{\eta^1 \ldots \eta^n}^{\xi} 
(\mathbb{F}_{p^l})$
is equal to 
$\mathbf{p}_{\eta^1 \ldots \eta^n}^{\xi} (p^l)$
for any prime number $p$ and any 
positive integer $l$.
\item\label{MultDegree}
The polynomial 
$\mathbf{p}_{\eta^1 \ldots \eta^n}^{\xi}$
is identically equal to zero or has degree
\begin{equation}\nonumber
\deg \mathbf{p}_{\eta^1 \ldots \eta^n}^{\xi} =
g_{| \eta^1 | \ldots | \eta^n |}^{| \xi |} +
\tfrac{1}{2} \bigl(
x_{\xi} +
x_{\eta^1} + \ldots + x_{\eta^n} 
\bigr) \; .
\end{equation}
\item\label{MultLR}
The leading coefficient of
the polynomial
$\mathbf{p}_{\eta^1 \ldots \eta^n}^{\xi}$ is
equal to
$d_{\eta^1 \ldots \eta^n}^{\xi}$.
\end{enumerate}
\end{theorem}
\begin{proof}
Statement \ref{MultDim} is proven in
\cite{Nakajima2001p} and 
\cite{Malkin2000b}. In these papers
the base field $k$ is assumed to be
$\mathbb{C}$, but the proofs work for 
arbitrary algebraically closed field.
\ref{MultIrr} is proven
in \cite{Malkin2000b} (it is also
a corollary of
\cite[Lemma 4.4]{Nakajima2001p} or 
\cite[Theorem 5.3]{VaragnoloVasserot2001}).
Statements \ref{MultDegree} and
\ref{MultLR} follow from
\ref{MultDim},
\ref{MultIrr}, and 
\ref{MultPolynomial}.

To prove \ref{MultPolynomial} one can
use the inductive argument in 
\cite[6.6]{Lusztig2000a} 
replacing the variety
(notation of \cite{Lusztig2000a})
$\Lambda^{s, \ast s}_{\mathbf{D},\mathbf{V}}$
(resp. $\Lambda^{s}_{\mathbf{D},\mathbf{V}}$) 
by
$\mathfrak{P}_{\mu^1 \ldots \mu^n}^{\lambda}$
(resp. a tensor product variety --
cf. \cite{Malkin2000b}).
To prove that the number of 
$\mathbb{F}_{p^l}$-rational
points of a tensor product variety
depends polynomially on $p^l$ note 
(cf. \cite{Malkin2000b}) that
a tensor product variety is a finite union 
of disjoint subsets each of which
is a vector bundle over a fibration with the 
base equal to a $\mathbb{Z} [I]$-graded 
partial flag variety and fibers isomorphic to
$\Lambda^{s}_{\mathbf{D}^1, 
\mathbf{V}^1; \mathbf{U}^1}
\times \ldots \times 
\Lambda^{s}_{\mathbf{D}^n, 
\mathbf{V}^n; \mathbf{U}^n}$ 
for some
$\mathbf{D}^1, \ldots , \mathbf{D}^n$,
$\mathbf{V}^1, \ldots , \mathbf{V}^n$,
$\mathbf{U}^1, \ldots , \mathbf{U}^n$
(again notation of \cite{Lusztig2000a}).
Thus it follows from 
\cite[Section 6]{Lusztig2000a} 
that the number of $\mathbb{F}_{p^l}$-rational
points of a tensor product variety
is given by a computable polynomial in $p^l$.
\end{proof}

\begin{remark}
The Hall polynomials
$\mathbf{h}_{\mu^1 \mu^2}^{\lambda}$
(cf. Section \ref{HallTheorem}) give 
structure constants of the Hecke
algebra of $GL (N)$ in the basis of
characteristic functions of the
double cosets.
The author does not know what connection
(if any) exists between the polynomials 
$\mathbf{p}_{\eta^1 \eta^2}^{\xi}$
and the Hecke algebra (of an extension of 
the algebraic group 
associated to $\mathfrak{g}$ by its Cartan
subgroup).
\end{remark}

\subsection{$A_{N-1}$ case}

In this section it is explained why
Theorem \ref{HallTheorem} is a
special case of Theorem
\ref{MultTheorem}. The argument
is quite standard 
(cf. \cite{KraftProcesi1979}).
Throughout the section $\mathfrak{g}$
is assumed to be $\mathfrak{sl}_{\sss N}$,
and the set of vertices $I$ 
of its Dynkin graph is identified with
$\{ 1, 2 , \ldots , N-1 \}$ in such a way that
two vertices are connected by an edge if 
and only if they are consequent integers.

Note that in the definition
\eqref{HallTensor} of 
$c_{\mu^1 \ldots \mu^n}^{\lambda}$
the number of non-zero parts 
of each of the partitions $\lambda ,
\mu^1 , \ldots , \mu^n ,$ is
assumed to be
less than $(N+1)$. Thus
representations of $k [[t]]$
involved in the definition of
$\mathfrak{N}^{\lambda}_{\mu^1 
\ldots \mu^n} (k)$
factor through 
$k [[t]]/k^Nk[[t]]$. 

Recall that $[1] \in \mathcal{F} (k)$ 
denotes the path of length zero beginning
and ending at the vertex $1$.

\begin{proposition}
Let $\mathcal{I}$ be a two-sided
ideal of $\Tilde{\mathcal{F}} (k)$ 
generated 
by $u_2 , u_3 , \ldots , u_{N-1}$.
\begin{enumerate}
\item\label{ANGenerator}
The factor algebra 
$\Tilde{\mathcal{F}} (k) / \mathcal{I}$
is generated as a $k$-algebra
by $(u_1 + \mathcal{I})$ and 
$([1] + \mathcal{I})$.
\item\label{ANRelation}
$\underbrace{[1] \circ [1] \circ \ldots 
\circ [1]}_{N} \in \mathcal{I}$ .
\item\label{ANIsomorphism}
A $k$-algebra homomorphism 
$\; \varphi : \;
k[[t]]/t^N k[[t]] \rightarrow 
\Tilde{\mathcal{F}} (k) / \mathcal{I}$
uniquely defined by 
$\varphi (1 + t^N k[[t]]) =
u_1 + \mathcal{I}$ and
$\varphi (t + t^N k[[t]]) =
[1] + \mathcal{I}$
is an isomorphism.
\end{enumerate}
\end{proposition} 

\begin{proof}
\ref{ANGenerator} and \ref{ANRelation} are
straightforward. Due to \ref{ANRelation} 
the map $\varphi$ is well-defined, and
due to \ref{ANGenerator} it is surjective. 
According to \cite{KraftProcesi1979}
there exists a faithful representation of
$k[[t]]/t^N k[[t]]$ which factors through 
$\varphi$. Thus $\varphi$ is injective. 
\ref{ANIsomorphism} follows. 
\end{proof}

Note that a representation of 
$\Tilde{\mathcal{F}} (k)$
factors through 
$\Tilde{\mathcal{F}} (k)/ \mathcal{I}$
if and only if its type
$\eta$
satisfies the condition:
$| \eta |_i = 0$
for $i=2, \ldots , N-1$.
Given such $\eta$ let
$( \mathbf{u}, \mathbf{v}) =
\varkappa^{-1} ( \eta )$ and let
$\nu (\eta )$ be an $N$-tuple
of integers given by  
$\nu (\eta )_1 = | \eta |_1 -
\mathbf{v}_1$,
$\nu (\eta )_i = \mathbf{v}_{i-1} -
\mathbf{v}_{i}$ for
$i = 2, \ldots , N-1$, and
$\nu (\eta )_N = \mathbf{v}_{N-1}$.
One has $\eta \in 
\mathcal{Q}^{++}_{\mathfrak{g}'}$
if and only if $\nu (\eta )$ 
is an ordered partition
(i.e. 
$0 \leq \nu (\eta )_i \leq \nu (\eta )_j$
for any $i \geq j$). Thus
$\nu$ provides a bijection
between types of representations of
$\Tilde{\mathcal{F}} (k) / \mathcal{I}$
and those of $k[[t]]/t^Nk[[t]]$. 
Moreover if a representation $\pi$ of
$\Tilde{\mathcal{F}} (k) / \mathcal{I}$
is of type $\eta$ then the 
representation $\pi \varphi$ of
$k[[t]]/t^Nk[[t]]$ is of type 
$\nu ( \eta )$, and 
one has the following 
bijection (isomorphism
of varieties if $k$ is
algebraically closed):
\begin{equation}\nonumber
\mathfrak{P}_{\eta^1 \ldots \eta^n}^{\xi} (k) 
\simeq
\mathfrak{N}_{\nu (\eta^1) \ldots 
\nu (\eta^n)}^{\nu (\xi )} (k) \; .
\end{equation}
Therefore Theorem \ref{HallTheorem}
follows from Theorem \ref{MultTheorem}.

\end{document}